\def\TheMagstep{\magstep1}	% Normal magnification
		% Changed to \magstep0 by \DoublepageOutput{TRUE} 
\def\PaperSize{letter}		% \PaperSize is used to

\magnification=\magstep1
%% A PlainTeX macro package

%%  MATH MACROS  %%
\let\:=\colon  
   \let\?=\overline

\let\Sum=\sum \def\sum{\Sum\nolimits}

\def\IC{{\bf C}} 
\def\IP{{\bf P}}

\def\pd #1#2{\partial#1/\partial#2}
\def\Pd #1#2{{\partial#1\over\partial#2}}
\def\sdep #1{#1^\dagger}	
\def\and{\hbox{ and }}
\def\Wf{\hbox{\rm W$_f$}}	\def\Af{\hbox{\rm A$_f$}}
\def\WF{\hbox{\rm W$_F$}}	\def\AF{\hbox{\rm A$_F$}}

\def\DONE{*!*}
\def\NextDef #1 {\def\NextOne{#1}%
 \ifx\NextOne\DONE\let\next\relax
 \else\expandafter\xdef\csname#1\endcsname{\TheOp}
  \let\next\NextDef
 \fi \next}
\def\TheOp{\mathop{\rm\NextOne}}
 \NextDef 
  Projan Supp Proj Sym Spec Hom cod Ker dist
 *!*
\def\TheOp{{\cal\NextOne}}
\NextDef 
  E F G H I J M N O R S
 *!*
\def\TheOp{\hbox{\rm\NextOne}}
\NextDef 
 A ICIS 
 *!*
 
 %%  STYLE MACROS  %%
%% Redefine \item to give greater indentation than AMSTeX
%   and the roman font within parentheses.   
\def\item#1 {\par\indent\indent\indent\indent \hangindent4\parindent
 \llap{\rm (#1)\enspace}\ignorespaces}
%% Define a similar macro without the hanging indentation for assertions
%% and that starts each part with an ordinary \parindent
 \def\inpart#1 {{\rm (#1)\enspace}\ignorespaces}
 \def\part {\par\inpart}
%% For displaying statements in italics with narrower margins and numbers
%\def\state{\smallskip\begingroup\narrower\noindent\it}

\catcode`\@=11		% make @ a letter temporarily

%% Modification of the PLAIN footnote macro for 8pt
\def\vfootnote#1{\insert\footins\bgroup
 \eightpoint %% only change
 \interlinepenalty\interfootnotelinepenalty
  \splittopskip\ht\strutbox % top baseline for broken footnotes
  \splitmaxdepth\dp\strutbox \floatingpenalty\@MM
  \leftskip\z@skip \rightskip\z@skip \spaceskip\z@skip \xspaceskip\z@skip
  \textindent{#1}\footstrut\futurelet\next\fo@t}

%%  ``Ties'' with a \thinspace for page numbers
\def\p.{p.\penalty\@M \thinspace}
\def\pp.{pp.\penalty\@M \thinspace}
%% For Roman parenthetical material in nonRoman text
\def\(#1){{\rm(#1)}}\let\leftp=(
\def\activeleftp{\catcode`\(=\active}
{\activeleftp\gdef({\ifmmode\let\next=\leftp \else\let\next=\(\fi\next}}

%% Sectioning
\def\sct#1\par
  {\removelastskip\vskip0pt plus2\normalbaselineskip \penalty-250 
  \vskip0pt plus-2\normalbaselineskip \bigskip
  \centerline{\smc #1}\medskip}

\newcount\sctno \sctno=0
\def\sctn{\advance\sctno by 1 
% {\bf\hbox to \parindent{\number\sctno.\hfil}#1.}\enspace\ignorespaces}
 \sct\number\sctno.\quad\ignorespaces}

%% Display numbers
\def\dno#1${\eqno\hbox{\rm(\number\sctno.#1)}$}
\def\Cs#1){\unskip~{\rm(\number\sctno.#1)}}

%% For setting results
\def\proclaim#1 #2 {\medbreak
  {\bf#1 (\number\sctno.#2)}\enspace\bgroup\activeleftp
\it}
\def\endproclaim{\par\egroup\medskip}
\def\pf{\endproclaim{\bf Proof.}\enspace}
 \def\prp{\proclaim Proposition }
\def\cor{\proclaim Corollary }	\def\thm{\proclaim Theorem }
\def\rmk#1 {\medbreak {\bf Remark (\number\sctno.#1)}\enspace}
\def\eg#1 {\medbreak {\bf Example (\number\sctno.#1)}\enspace}

%% PAGE LAYOUT
\parskip=0pt plus 1.75pt \parindent10pt
\hsize29pc
\vsize44pc
\abovedisplayskip6pt plus6pt minus2pt
\belowdisplayskip6pt plus6pt minus3pt

\def\TRUE{TRUE}	% For Boolean tests
\ifx\DoublepageOutput\TRUE \def\TheMagstep{\magstep0} \fi
\mag=\TheMagstep

% CENTER TEXT ON PAGE
	% additional vertical adjustment
\newskip\vadjustskip \vadjustskip=0.5\normalbaselineskip
\def\centertext
 {\hoffset=\pgwidth \advance\hoffset-\hsize
  \advance\hoffset-2truein \divide\hoffset by 2\relax
  \voffset=\pgheight \advance\voffset-\vsize
  \advance\voffset-2truein \divide\voffset by 2\relax
  \advance\voffset\vadjustskip
 }
\newdimen\pgwidth\newdimen\pgheight
\def\letter{letter}\def\AFour{AFour}
\ifx\PaperSize\letter
 \pgwidth=8.5truein \pgheight=11truein 
 \message{- Got a paper size of letter.  }\centertext 
\fi
\ifx\PaperSize\AFour
 \pgwidth=210truemm \pgheight=297truemm 
 \message{- Got a paper size of AFour.  }\centertext
\fi

%% TWO-COLUMN LANDSCAPE FORMAT
% Modified from the TeX book, p. 257.
 \newdimen\fullhsize \newbox\leftcolumn
 \def\fulline{\hbox to \fullhsize}
\def\doublepageoutput
{\let\lr=L
 \output={\if L\lr
          \global\setbox\leftcolumn=\columnbox \global\let\lr=R%
        \else \doubleformat \global\let\lr=L\fi
        \ifnum\outputpenalty>-20000 \else\dosupereject\fi}%
 \def\doubleformat{\shipout\vbox{%
        \fulline{\hfil\hfil\box\leftcolumn\hfil\columnbox\hfil\hfil}%
				}%
		  }%
 \def\columnbox{\vbox
   {\makeheadline\pagebody\makefootline\advancepageno}%
   }%
 \fullhsize=\pgheight \hoffset=-1truein
 \voffset=\pgwidth \advance\voffset-\vsize
  \advance\voffset-2truein \divide\voffset by 2
  \advance\voffset\vadjustskip
 \let\firstheadline=\hfil
 
% \null\vfill\nopagenumbers\eject\pageno=1\relax % to put page on right
}
\ifx\DoublepageOutput\TRUE \doublepageoutput \fi

%% ADDITIONAL FONTS
 \font\twelvebf=cmbx12		% For title
 \font\smc=cmcsc10		% For authors' names

%% EIGHT POINT TYPE FOR FOOTNOTES AND REFERENCES
\def\eightpoint{\eightpointfonts
 \setbox\strutbox\hbox{\vrule height7\p@ depth2\p@ width\z@}%
 \eightpointparameters\eightpointfamilies
 \normalbaselines\rm
 }
\def\eightpointparameters{%
 \normalbaselineskip9\p@
 \abovedisplayskip9\p@ plus2.4\p@ minus6.2\p@
 \belowdisplayskip9\p@ plus2.4\p@ minus6.2\p@
 \abovedisplayshortskip\z@ plus2.4\p@
 \belowdisplayshortskip5.6\p@ plus2.4\p@ minus3.2\p@
 }
\newfam\smcfam
\def\eightpointfonts{%
 \font\eightrm=cmr8 \font\sixrm=cmr6
 \font\eightbf=cmbx8 \font\sixbf=cmbx6
 \font\eightit=cmti8 
 \font\eightsmc=cmcsc8
 \font\eighti=cmmi8 \font\sixi=cmmi6
 \font\eightsy=cmsy8 \font\sixsy=cmsy6
 \font\eightsl=cmsl8 \font\eighttt=cmtt8}
\def\eightpointfamilies{%
 \textfont\z@\eightrm \scriptfont\z@\sixrm  \scriptscriptfont\z@\fiverm
 \textfont\@ne\eighti \scriptfont\@ne\sixi  \scriptscriptfont\@ne\fivei
 \textfont\tw@\eightsy \scriptfont\tw@\sixsy \scriptscriptfont\tw@\fivesy
 \textfont\thr@@\tenex \scriptfont\thr@@\tenex\scriptscriptfont\thr@@\tenex
 \textfont\itfam\eightit	\def\it{\fam\itfam\eightit}%
 \textfont\slfam\eightsl	\def\sl{\fam\slfam\eightsl}%
 \textfont\ttfam\eighttt	\def\tt{\fam\ttfam\eighttt}%
 \textfont\smcfam\eightsmc	\def\smc{\fam\smcfam\eightsmc}%
 \textfont\bffam\eightbf \scriptfont\bffam\sixbf
   \scriptscriptfont\bffam\fivebf	\def\bf{\fam\bffam\eightbf}%
 \def\rm{\fam0\eightrm}%
% \tt \ttglue=0.5em plus0.25em minus0.15em
 }

%% HEADLINE STYLE
\def\today{\ifcase\month\or	% From the TeX book p. 406
 January\or February\or March\or April\or May\or June\or
 July\or August\or September\or October\or November\or December\fi
 \space\number\day, \number\year}
\nopagenumbers
\headline={%
  \ifnum\pageno=1\firstheadline
  \else
    \ifodd\pageno\oddheadline
    \else\evenheadline\fi
  \fi
}
\let\firstheadline\hfill
\def\oddheadline{\eightpoint \rlap{\today}
 \hfil\headtitle\hfil\llap{\folio}}
\def\evenheadline{\eightpoint\rlap{\folio}
 \hfil\author\hfil\llap{\today}}
\def\headtitle{\title}

%% REFERENCING
	% to introduce the keys in order
 \newcount\refno \refno=0	 \def\NoKey{*!*}
 \def\MakeKey{\advance\refno by 1 \expandafter\xdef
  \csname\TheKey\endcsname{{\number\refno}}\NextKey}
 \def\NextKey#1 {\def\TheKey{#1}\ifx\TheKey\NoKey\let\next\relax
  \else\let\next\MakeKey \fi \next}
 \def\RefKeys #1\endRefKeys{\expandafter\NextKey #1 *!* }
	% to set references
\def\SetRef#1 #2,#3\par{%
 \hang\llap{[\csname#1\endcsname]\enspace}%
  \ignorespaces{\smc #2,}
  \ignorespaces#3\unskip.\endgraf
 }
 \newbox\keybox \setbox\keybox=\hbox{[8]\enspace}
 \newdimen\keyindent \keyindent=\wd\keybox
\def\references{%\vskip-\smallskipamount
  \bgroup   \frenchspacing   \eightpoint
   \parindent=\keyindent  \parskip=\smallskipamount
   \everypar={\SetRef}}
\def\endreferences{\egroup}

%% SERIALS
 \def\serial#1#2{\expandafter\def\csname#1\endcsname ##1 ##2 ##3
  {\unskip\ #2 {\bf##1} (##2), ##3}}
 \serial{ajm}{Amer. J. Math.}
  \serial {aif} {Ann. Inst. Fourier}
 \serial{asens}{Ann. Scient. \'Ec. Norm. Sup.}
 \serial{comp}{Compositio Math.}
 \serial{conm}{Contemp. Math.}
 \serial{crasp}{C. R. Acad. Sci. Paris}
 \serial{dlnpam}{Dekker Lecture Notes in Pure and Applied Math.}
 \serial{faa}{Funct. Anal. Appl.}
 \serial{invent}{Invent. Math.}
 \serial{ma}{Math. Ann.}
 \serial{mpcps}{Math. Proc. Camb. Phil. Soc.}
 \serial{ja}{J. Algebra}
 \serial{splm}{Springer Lecture Notes in Math.}
 \serial{tams}{Trans. Amer. Math. Soc.}

	% modified \cite code from AMSTeX
\def\UThin{\penalty\@M \thinspace\ignorespaces}
	% unbreakable \thinspace for use after periods
\def\relaxnext@{\let\next\relax}
\def\cite#1{\relaxnext@
 \def\nextiii@##1,##2\end@{\unskip\space{\rm[\SetKey{##1},\let~=\UThin##2]}}%
 \in@,{#1}\ifin@\def\next{\nextiii@#1\end@}\else
 \def\next{{\rm[\SetKey{#1}]}}\fi\next}
\newif\ifin@
\def\in@#1#2{\def\in@@##1#1##2##3\in@@
 {\ifx\in@##2\in@false\else\in@true\fi}%
 \in@@#2#1\in@\in@@}
\def\SetKey#1{{\bf\csname#1\endcsname}}

\catcode`\@=12 %\active  %at signs are no longer letters

\def\title{ Nilpotents, Integral Closure and Equisingularity conditions }
\def\author{Terence Gaffney}
\RefKeys  G0 G1 G2 G3 G4 GK GK2 KT1 LT
 \endRefKeys

\def\topstuff{\leavevmode
 \bigskip\bigskip
 \centerline{\twelvebf \title}
 \bigskip
 \centerline{\author}
 \medskip\centerline{\today}

\bigskip\bigskip}

\topstuff

\sct Introduction

The motivation for this note comes from the fact that if you want an equisingularity theory for general singularities
parallel to the complete
 intersection case, you must admit nilpotents into your structure sheaves.

 Here is a simple example. Consider a family whose members 
are 2-planes in 4-space meeting transversally at a single point,
 with one plane fixed and the other moving.  In example 1 we show that if we take a hyperplane section
 of a member of this family we get a non-reduced space; if we take a hyperplane section
 of the total space of this family, which contains the parameter space, we again get a space with nilpotents.

 We could take the reduced structure on the intersection, but the fat point says something interesting.
 It exists because of the possibilty of pulling the two lines which make up the hyperplane 
section apart by a flat deformation. So the question arises, how do we have to change the 
theory to take nilpotents into account? In \cite {G1}, we gave an algebraic condition 
for the Whitney conditions to hold in terms of the integral closure of the Jacobian module associated to the defining 
equations of the analytic set. Assume $X$ is equidimensional as a set. In this note, in section 1 we show that 
we don't need to take the reduced struture for $X$ in checking this algebraic criterion, we can use the generators for any ideal 
which will give the reduced structure on $X$ at smooth points of $X$.  

This is useful in applications where often the set of equations constructed by some geometric 
process are only known to give the reduced structure at the set of smooth points.

So the integral closure approach is blind to nilpotents just as the Whitney conditions are. As a consequence, in section 2, in Theorem 2.2, we can prove easily 
that the Whitney conditions are preserved by intersection with a generic hyperplane using the theory of integral closure, and we can describe precisely when a hyperplane is generic in this sense.

In the course of studying what it means for a hyperplane to be generic for this situation, we are led to prove a result in the case of families of complete intersections with isolated singularities, relating the limiting tangent hyperplanes to a complex analytic set $X$ and the fibers of $X$ over $Y$ (Theorem 2.6). This allows an improvement of Theorem 2.2 in this case.

The author thanks Steven Kleiman for helpful conversations. This project grew out of a collaboration with David Trotman and Leslie Wilson; the author thanks both of his collaborators for their interest and support.

\sctn Equisingularity and Nilpotents

We start with a motivating example.

\eg 1 Let $X$ be the germ of an analytic set defined at the origin in $\IC^5$  with reduced structure by 
$$\{x(w+ty)=0, xz=0, y(w+ty)=0, yz=0\}.$$

The structure is reduced because the collection $\{ x,y,z, w+ty\}$ is an R-sequence in $\O_5$,
 the local ring of germs of analytic functions at the origin
in $\O_5$, so the only relations between these elements are the Kozul relations. The projection from $\IC^5$ to the $t$-axis, makes $X$ into a family of sets; the fiber over 
$0$ , $X(0)$ is defined by 
$$\{xw=0, xz=0, yw=0, yz=0\},$$
\noindent and is also reduced for the  same reason. So, $X$ is a family consisting of 2 two-planes, one plane fixed and the other moving.
The $t$-axis is the parameter space of the family
 and is also the singular set of 
$X$.

 However, if we intersect with a hyperplane of the form $x=ay+bz+cw$,  $c\ne 0$, 
the simplest structure to use on the intersections is that given by adding the equation of $H$ to the equations defining
$X$ and $X(0)$. If we use the equation for $H$ to eliminate $x$, then this structure on $X\cap H$ contains the nilpotent $(b/c)z+w+ty$, while the structure
on $X(0)\cap H$ contains the nilpotent $(b/c)z+w$. (If we eliminate $x$, the "simple structure`` is defined by the product of the ideals
$I_1=(bz+cw,y)$ and $I_2=(w+ty,z)$, while the reduced structure is defined by the intersection of the two ideals, and the intersection contains the element
$(b/c)z+w+ty$ which is not in the product.)

As we shall see, many equisingularity conditions are described using information from the Jacobian module of $X$. If $X,x\subset \IC^n$ is defined as
$F^{-1}(0)$ where $F:\IC^n,x\to \IC^p,0$, then the Jacobian module of $X$ with the structure defined by $F$, is just the submodule of $\O^p_{X,x}$ generated by the partial derivatives of $F$ and denoted $JM(F)$. Sometimes
we have $X,x\subset \IC^k\times \IC^n$, with $y$ a set of coordinates on $\IC^k$ and $z$ a set of coordinates on  $\IC^n$. In this case, we let $JM_z(F))$ denote the submodule of $\O^p_{X,x}$ generated by the partial derivatives of $F$ with respect to $z$, $JM_y(F)$ denotes the submodule of $\O^p_{X,x}$ generated by the partial derivatives of $F$ with respect to $y$. The example of the previous paragraph shows that if we wish to pass information from the Jacobian module of $X$ to that of $X\cap H$ easily, 
we must deal with nilpotent structures on both the total space and the members of the family.

We now introduce the notion of integral closure, and begin to consider the connection between nilpotents and this idea.

Recall that if $R$ is a commutative, unitary ring, $I$ an ideal in $R$, $f$ an element of $R$, then $f$ is integrally dependent
on $I$ if there exists a positive integer $k$ and elements $a_j$ in
$I^j,$ so that $f$ satisfies the relation
$f^k+a_1f^{k-1}+\dots+a_{k-1}f+a_k=0$ in
$R.$

The integral closure of an ideal $I$, denoted $\bar I$, consists of all elements integrally dependent on it.

\prp 2 Supose $I$ is an ideal in a commutative, unitary ring $R$. If $h$ is a nilpotent element of $R$, then $h\in
\bar I$.

\pf Since $h$ is nilpotent it satisfies an equation of the form $T^k=0$ for some $k$.

\prp 3 Suppose $R$ as above, and $R_r$ is the ring gotten 
by modding out by nilpotents, p the projection map. Suppose $I$ is an ideal in $R$, $h$ an element of $R$, then
$\bar I$ contains $h$ iff $\overline{(p(I))}$
 contains $ p(h)$.

\pf If $\bar I$ contains $h$, then we just reduce mod nilpotents and the result follows. Suppose $\overline {(p(I))}$
 contains $ p(h)$. Then $ p(h)$ satisfies an equation of integral dependence. Choosing representatives of the coefficients of the equation
we have $h$ satisfies an equation of the form $P(h)=g$ where $g$ is nilpotent. Then consider $P^k(h)=g^k$. For $k$ large enough,
$g^k=0$, and $P^k$ is a polynomial of the desired form, because it is a product of polynomials of the desired form.

We can also talk about the integral closure of a module.  Suppose $M$ is a submodule of $E=R^k$, $k$ at
least 1.   Let $\rho: E\to SE$ be the inclusion of $E$ into its symmetric
algebra; then $h\in E$ is integrally dependent on $M$ if $\rho(h)$ is integrally dependent on the ideal generated by $\rho(M)$. 

If $N\subset M$ are a pair of modules such that $\bar N=\bar M$, then we say that $N$ is a reduction of $M$.

\cor 4 Suppose $R$ is a commutative, unitary ring,  $R_r$ is the ring gotten 
by modding out by nilpotents, p the projection map from $R^k$ to $R^k_r$,  suppose $M$ is a submodule in $R^k$,
 $h$ an element of $R^k$, then
$\bar M$ contains $h$ iff $\overline {(p(M))}$
 contains $ p(h)$.

\pf Translate from modules to ideals and use proposition 1.2.

 There is a useful criterion (curve criterion) for checking if an element is in the integral closure of a module, which also allows us to define the notion of strict dependence. Suppose $X,x$ is a complex analytic germ, $M$ a submodule of
$\O^p_{X,x}$. Then $h\O^p_{X,x}$ is in the integral closure
of $M$ (resp. strictly dependent on $M$) iff for all
$\phi :{\IC},0 \rightarrow X,x$, $ h \circ \phi \in (\phi^{*}M){\O}_{1}$ (resp. if for all $\phi : {\IC},0
\rightarrow X,x$ we have $h \circ\phi\in m_{1}\phi^{*}M$, where $m_{1}$ is the
maximal ideal in ${\O}_{1}$). (For a discussion of the curve criterion cf. \cite{G1}.)

We denote the set of elements strictly
dependent on $M$ by ${M}^\dagger$. \par\vskip .2in

There is also a version of Nakayama's lemma for integral closure based on the curve criterion.

\prp 5 Suppose $N\subset M\subset \O^p_{X,x}$, and $\?{N+M^\dagger}=\?M$, then $N$ is a reduction of $M$.

\pf  Use the curve criterion and Nakayama's lemma.

The theorem of Rees is known to hold if even if the local ring is not reduced. We have:

 \prp 6 Suppose $X$ is the germ of an analytic space, such that $X_{red}$ is equidimensional, 
$I$ and $J$ two ideals in the local ring of $X$. 
Suppose $e(I)=e(J)$,  $I\subset J$, then $J$ in $\bar I$.

\pf Cf \cite {KT1}

In fact, the multiplicity of $I$ is the same
 as the multiplicity of $(p(I))$, when the structure on $R$ is generically reduced. 
Here is an argument by Steven Kleiman proving this fact.

We have the map of Rees algebras
$R(I)\to R(p(I))$ is surjective and has a nilpotent kernel, so the blowup of
Spec$(R_r)$ is simply the reduction of the blowup of Spec$(R)$, and the
exceptional divisor of the first induces to that of the second; the
equation now follows from the projection formula applied to the
inclusion of the reduced scheme in the nonreduced one.

For the various stratification conditions we still need to do some work, because these depend on the Jacobian module, which means we must 
differentiate our functions which are giving us the non-reduced structure. 

The setup for our theorems is as follows: $X,0\subset \IC^{k+n},0$ is the germ of a complex analytic set at the origin, $Y=\IC^k\times 0\subset X $, $(y_1,\dots, y_k)$ coordinates on $\IC^k$, $(z_1,\dots, z_n)$ coordinates on $\IC^n$. $X_0$ denotes the smooth points of $X$, $X$ with reduced structure. We denote the ideal sheaf generated by the  $(z_1,\dots, z_n)$  by $m_Y$. We assume $X_{red}$ defined by an ideal
$J$, $I$ an ideal such that $V(I)=X$, and $I$ gives the reduced structure on $X$ off the singular set of $X$. Let  $f$ and $g$ be map germs whose components are a set of generators for $J$ and $I$. Let $F:X,Y\to\IC,0$.

\thm 7 Suppose $X,0\subset \IC^{n+k},0$ is the germ of an analytic space, such that $X_{red}$ is equidimensional, $X$,  $I$ and $J$, $Y$, $f$, $g$ as in the set-up above then:
$$\Pd f {y_i}\in \?{JM_z(f)} \hskip 1em{\rm {iff}}\hskip 1em \Pd g {y_i}\in \?{JM_z(g)}.$$ 

\pf We will see shortly that the implication $\Pd f {y_i}\in \?{JM_z(f)}$ implies $\Pd g {y_i}\in \?{JM_z(g)}$ is trivial, 
so we will concentrate on the other implication.  We know that the inclusion $\Pd g {y_i}\in \?{JM_z(g)}$ holds with the reduced structure,  
and we will show the desired implication using the curve criterion. We need only use curves whose 
image (except for the origin) lies in the smooth part of $X$. This is because  inclusion at the module level is equivalent
to inclusion of the corresponding Fitting ideals, and to check this we only need the kind of curves we are considering.

Since $I\subset J$, we can form the quotient sheaf; this is supported on the singular set of $X$ by hypothesis, so 
$I\supset K^rJ$ for some $r$ where $K$ defines the singular locus of X, by the Nullstellensatz.

Suppose the number of components of $f$ is $p$ while that of $g$ is $p'$. Then there is a $p'\times p$ matrix $H_0$ such
that $H_0(f)=(g)$. Since $I\supset K^rJ$, there exist matrices $H_1$, $H_2$ such that $H_1(g)=H_2(f)$, where the entries of
$H_2(f)$ are $q_{(j-1)p+t}=k^r_{j}f_{t}$, $k_j$ a set of generators of $K$, $1\le t\le p$. The matrix $H_2$, whose entries are 
$h_{(j-1)p+t,t}=k^r_{j}$, $1\le t\le p$ has maximal rank at all smooth points of $X$.

If we differentiate each of our column vectors we obtain:
$$Dg=H_0 Df$$

$$H_1Dg=H_2Df$$
working over the local ring of $X_{red}$. Now we pull back by a curve $\phi$. Since the rest of the argument 
is independent of the dimension of $Y$
we assume dim $Y=1$ for notational simplicity.

Now the inclusion $\Pd g {y}\in \?{JM_z(g)}$ implies that $\Pd g {y}\circ\phi\in {JM_z(g)}\circ \phi$. This is equivalent to the existence of $(1, v(t))$ such that
$$(Dg\circ \phi) (1, v)=0.$$

It's now clear that $\Pd f {y}\in \?{JM_z(f)}$ implies $\Pd g {y}\in \?{JM_z(g)}$ 
 as $0=(Df\circ \phi) (1, v)$ implies $0=H_0\circ\phi (Df\circ \phi) (1, v)=(Dg\circ \phi )(1, v)=0$.

Suppose $(Dg\circ \phi )(1, v)=0$. Then $H_1\circ\phi (Dg\circ \phi) (1, v)=H_2\circ\phi (Df\circ\phi)(1,v)$. Since $H_2\circ\phi$
is invertible for $t\ne 0$, it follows that $0=(Df\circ\phi)(1,v)$.

It is now easy to see how to improve our theorems relating integral closure and the various stratification conditions.

\thm 8 Suppose $X,0\subset \IC^{n+k},0$ is the germ of an analytic space, such that $X_{red}$ is equidimensional, $X$,  $I$ and $J$, $Y$, $f$, $g$, $F$ as in the set-up above then:

$$1) \Pd f {y_i}\in \?{m_YJM_z(f)} \hskip 2pt\forall i \hskip 2pt{\rm {iff}}\hskip 2pt \Pd g {y_i}\in \?{m_YJM_z(g)} \hskip 2pt \forall i\hskip 2pt
 {\rm iff}  \hskip 2pt  {\rm the \hskip 2pt pair }\hskip 2pt(X_0, Y)\hskip 2pt {\rm satisfies} \hskip 1pt W.$$

$$2) \Pd {(f ,F)}{y_i}\in \sdep{JM_z(f,F)}\hskip 2pt \forall i \hskip 2pt{\rm {iff}}\hskip 2pt \Pd {(g ,F)}{y_i}\in \sdep{JM_z(g,F)}\hskip 2pt \forall i \hskip 2pt$$

$$\hskip 2in {\rm iff}  \hskip 2pt  {\rm the \hskip 2pt pair }\hskip 2pt(X_0, Y)\hskip 2pt {\rm satisfies} \hskip 1pt \AF.$$

$$3) \Pd {(f ,F)}{y_i}\in \?{m_YJM_z(f,F)}  \hskip 2pt\forall i \hskip 2pt{\rm {iff}}\hskip 2pt \Pd {(g ,F)} {y_i}\in \?{m_YJM_z(g,F)}\hskip 2pt \forall i \hskip 2pt$$
$$\hskip 2in{\rm iff}  \hskip 2pt  {\rm the \hskip 2pt pair }\hskip 2pt(X_0, Y)\hskip 2pt {\rm satisfies} \hskip 1pt \WF.$$
\pf
Here is the argument for the first implication in 1).We know $Dg=H_0 Df$,  $H_1Dg=H_2Df$. By linearity we know
$z_i\Pd g {z_j}=H_0z_i\Pd f {z_j}$, so $$H_1H_0[\Pd f y, [z_i\Pd f {z_j}]]=H_2[\Pd f y, [z_i\Pd f {z_j}]],$$
 and now the argument proceeds as in Theorem 6. The equivalence of $ \Pd f {y_i}\in \?{m_YJM_z(f)} $ and W is Theorem 2.5 of \cite{G1} p309.

For the \Af case, define $\hat{H_0}$ to be the matrix such that the lower right corner entry is $1$, the other entries in the last row  and column 0, and the rest of the matrix
is $H_0$. Then $\hat{H_0}(f,F)=(g,F)$, and differentiating we get

$$\hat{H_0}D(f,F)=D( g,F).$$

Extend $H_2$ to $\hat{H_2}$, $H_1$ to $\hat{H_1}$  as we extended $H_0$ to $\hat{H_0}$. It is easy to see that

$$\hat{H_1}D(g,F)=\hat{H_2}D(f,F).$$

Now the argument goes as before. (Though we need strict dependence for \Af , this is not a problem, because we can assume
$v$ vanishes at the origin.) The \Wf argument combines the last two arguments. The equivalences with the stratification conditions comes from Lemma 5.1 p565 of \cite{GK} in the \AF \hskip 2pt case and from Proposition 2.1 p36 \cite{GK2}.

\sctn  Generic plane sections and the Whitney conditions

In this section, we apply the results on nilpotents to the study of sections of $(X_0, Y)$.

Setup: Given $X$, $Y$ as in the setup before Theorem 7, assume the pair $(X_0,Y)$ satisfies condition $W$. Consider a sequence of hyperplanes $H_1,H_2,\dots,H_n$, such that $P_i=\bigcap\limits_{j=1}^i H_j$ is a plane of codimension $i$ for all $i$, all $H_i\supset Y$, and $H_i$ defined by a linear form $F_i$. Note the set of hyperplanes in $Y\times\IC^n$ which contain $Y$ are parametrised by $\IP^{n-1}$.

We want to characterize those sequences for which  $(X_i=P_i\cap X_0, Y)$ satisfies $W$, and in which $(X_i=P_i\cap X_0, Y)$  is generic in a certain sense, which we will develop. 

To study the relation of the hyperplanes and the pair $(X_0,Y)$, we use a variant of the Grassman modification (\cite {G2}). In $\IC^k\times\IC^n\times \hat \IP^{n-1}$ consider the incidence variety $\widetilde{\IC^k\times\IC^n}=\{((y,z),H)| (y,z)\in H, H\in \IP^{n-1}\}$. Denote the projection to $\IC^k\times\IC^n$ by $\beta$, define  $\widetilde X$,  to be $\beta^{-1}(X)$. We can assume that we will be working with hyperplanes with equations $z_n=a_1z_1+\dots+ a_{n-1}z_{n-1}$. Then we can use $(y, z_1,\dots, z_{n-1}, a_1,\dots a_{n-1})$ as coordinates on  $\widetilde{\IC^k\times\IC^n}$, and in these coordinates $\beta(y,z,a)=(y,z, a_1z_1+\dots+ a_{n-1}z_{n-1})$. In these coordinates, we give  $\widetilde X$
 the scheme structure defined by $\beta^*(I)$ where $I$ defines $X$. Since $\beta$ is a submersion off $\IC^k\times 0\times  \hat \IP^{n-1}$, this structure will be reduced off of $Y\times  \hat \IP^{n-1}$, assuming $I$ gives the reduced structure off $Y$. (Otherwise it will be reduced off $\beta^{-1}S$, where $S$ is the set of points  where $I$ fails to give the reduced structure.)  It is  known that $(\widetilde X_0, Y\times  \hat \IP^{n-1)}$ satisfies W off some Z-open dense subset $V$ of  $Y\times  \hat \IP^{n-1}$. 
 Let $U\subset  \hat \IP^{n-1}=V\cap 0\times \hat \IP^{n-1}$. If $H$ is an element of $U$, it follows that  $(\widetilde X_0, Y\times  \hat \IP^{n-1)}$ satisfies W at $(0,H)$. We say such $H$ are W-generic hyperplanes for $(X_0,Y)$. A sequence of hyperplanes as in the setup is W-generic if each $H_i$ is W-generic for $(X_{i-1},Y)$.
 
 Since the curve criterion is so helpful in dealing with integral closure questions, it is helpful to know when a curve on $X,0$ lifts to $\widetilde X,0,H$, $H$ a hyperplane containing $Y$. Given a curve $\phi:\IC,0\to X,0$ the limiting $Y$-secant of $\phi$ is the limit as $t$ tends to $0$ of the line determined by $<\phi_{k+1}(t),\dots,\phi_{k+n}(t)>$.

\prp 1 Suppose $X,Y,0$ is the germ of a pair of complex analytic sets, $(X,Y)$ as in the setup before theorem 1.6. Then $\phi:\IC,0\to X,0$, where the image of $\phi$ does not lie in $Y$, lifts to $\widetilde X,0,H$, $H$ a hyperplane containing $Y$, if and only if $H\supset l$, where $l$ is the limiting $Y$-secant of $\phi$.

%Then $H$ is a limiting tangent hyperplane to $X$ at the origin, if and only if there exists $\phi :{\IC},0 \rightarrow X,x$, 
%such that at each point $\phi(t)$ there exists a tangent hyperplane to $X$, and $H$ is the limit of these tangent planes, while the limiting tangent line to $\phi$ at $x$ lies in $H$.

%\pf Suppose $H$ is a limiting tangent hyperplane  to $X,x$. Then we can find $\phi$, since $H$ is a point of the conormal space $C(X)$, such that  at each point $\phi(t)$ there exists a tangent hyperplane to $X$, and $H$ is the limit of these tangent planes. Each hyperplane contains the tangent line to $\phi$ at $\phi(t)$, so the limit contains the limiting line. The other implication is obvious.

\pf Suppose the equation of $H$ is $z_n=\sum_{i=1}^{n-1} a_iz_i$. Such a $\phi$ has an extension if and only if  there exists $a_i(t)$, $a_i(0)=a_i$, such that $\phi_{k+n}(t)=\sum_{i=1}^{n-1} a_i(t)\phi_{k+i}(t)$. It is easy to see that 
this last equation holds iff 

$$\mathop{\lim}\limits_{t\rightarrow 0} (1/t^{k})( \phi_{k+n}(t)-\sum_ia_{i}\phi_{k+i}(t)) =0$$
where  k is the minimum of the 
order of the first non-vanishing term in $\{\phi_{k+i}(t)\}$.

%Note that if $H$ was not a hyperplane, but merely contained in a limit of hyperplanes, there would be no reason why $H$ should contain the limiting tangent line.

We will give two characterizations of the W-generic hyperplanes for $(X_0 ,Y)$, one geometric, one algebraic. If $H$ is a hyperplane on $\IC^k\times\IC^n$ denote the submodule
of $JM(f)$ generated by the directional derivatives of $f$ in directions tangent to $H$ by $JM(f)_H$. If $H$ contains $Y$, then let $JM_z(f)_H$ denote the submodule of 
$JM_z(f)$ generated by directional derivatives of $f$ in directions tangent to $H\cap \IC^n$.

\thm 2 In the setup of this section, $H$ is W-generic for $(X_0, Y)$ iff either (hence both) of the following conditions hold

1) $H$ is not a limiting tangent hyperplane to $X$ at the origin.

2) $JM(f)_H$ is a reduction of  $JM(f)$  as $\O_{X}$ modules.

\pf  1) and 2) are equivalent by  \cite {G0}. Further, if $H$ is not a limiting tangent hyperplane, then it is clear that $(X\cap H)_0=X_0\cap H$. For $(X\cap H)_0\subset X_0\cap H$, while if there is a curve of singular points on $X_0\cap H$, then $H$ is a tangent hyperplane at such points, so $H$ is a limiting tangent hyperplane. 

The proof from here is similar to that of Theorem 2.9 of \cite{G2}.

Denote $f\circ \beta$ by $G$. A calculation using the chain rule shows that:

$$\pd G{a_{i}}=z_i\pd f{z_n}\circ\beta$$
$$\pd G{y_i}=\pd f{y_i}\circ\beta$$
$$JM_z(G)=(\pd f{z_j}\circ\beta+a_{j}\pd f{z_n}\circ\beta)$$
 By Theorem 1.8 1), $H$ is generic if and only if 
 $$\pd G{a_{i}}, \pd G{y_i}\in \? {\beta^*(m_Y)JM_z(G)}.$$
 
Using the curve criterion, and Nakayama's lemma,  we see that $z_i\pd f{z_n}\circ\beta\in \? {(\beta^*(m_Y)JM_z(G)}$, if and only if $\pd f{z_n}\circ\phi\in \phi^* JM_z(f)_H$, for $\phi$ any curve on $X_0$, whose limiting $Y$-secant is in $H$. (These are the curves on $X_0,0$ which lift to $\widetilde X_0, (0,H)$ by 2.1)

Suppose $H$ is not a limiting tangent hyperplane. Then $JM(f)_H$ is a reduction of $JM(f)$. Since $(X_0,Y)$ satisfies W at the origin, we have that 
$JM_y(f)\subset JM(f)^\dagger$ 

This implies that 
$$\?{JM_z(f)_H+ JM(f)^\dagger}=\?{JM(f)}$$
so, $JM_z(f)_H$ is a reduction of $JM(f)$ by 1.5. Then $JM_y(f)\subset \?{m_YJM_z(f)_H}$ and 
$\?{JM_z(f)_H}=\?{JM_z(f)}$ imply that  $z_i\pd f{z_n}\circ\beta\in \? {(\beta^*(m_Y)JM_z(G)}$ and $ \pd G{y_i}\in \? {\beta^*(m_Y)JM_z(G)}$, which shows $H$ is W-generic.

Suppose $H$ is W-generic. Then $JM_y(f), \pd f{z_n}\circ\phi\in \phi^* JM_z(f)_H$ holds for all $\phi$,whose limiting $Y$-secant is in $H$, again by the curve criterion and Nakayama's lemma. This checks the condition of 2) for all curves whose limiting $Y$-secant is in $H$. In particular $H$ cannot be a limit of tangent hyperplanes along such a curve. Now suppose $H$ is a limiting tangent hyperplane along some curve $\phi$; then since the pair $(X_0,Y)$ satisfy W and hence Whitney B, the limiting secant line of $\phi$ must be in the limiting tangent plane, which gives a contradiction.

\cor 3 Suppose $H$ is a W-generic hyperplane for $(X_0,Y)$, then $(X_0\cap H, Y)$ satisfies condition W.

\pf By the above proof $JM_z(f)_H$ is a reduction of $JM(f)$, this implies that  $JM_y(f)\?{m_YJM_z(f)_H}$. Let $G=(f, H)$; $G$ defines $X\cap H$ with possibly non-reduced structure.
Then the inclusion $JM_y(f)\subset \?{m_YJM_z(f)_H}$ continues to hold when we restrict to $\O^p_{X\cap H}$. This implies the inclusion $JM_y(G)\subset \?{JM_z(G}$ since $H$ is independent of $y$. The result follows from 1.8 1).

\rmk 4 The implication $H$ not a limit of tangent hyperplanes implies W-generic was proved by L\^e and Teissier (\cite {LT}), using the aureole and the theory of polar varieties. They also worked with the integral closures of Fitting ideals, as the theory of the integral closure of modules was not available then.

\cor 5 Suppose $(X_0,Y)$ satisfy W, $\{H_i\}$ are a sequence of hyperplanes as in the set up for this section, then $\{H_i\}$ is a W-generic sequence if and only if
each $H_i$ is not a limiting tangent plane at the origin for $X_{i-1}$.

\pf Apply Theorem 2 inductively to each  $X_{i-1}$.

Existing technology allows for a significant improvement if $X$ is a family of complete intersections with isolated singularities (ICIS). 

Given an equidimensional  complex analytic germ $X\subset C^k\times C^n,0$ containing $Y=C^k\times 0$, suppose that the pair $X_0,Y$ satisfies condition WA at the origin. This means that every limit of tangent hyperplanes at the origin contains $Y$.
 There is always a bijection $\iota$ which takes hyperplanes containing $Y$ to hyperplanes through the origin in $C^n$, given by intersecting these hyperplanes with $0\times C^n$. We can ask does this bijection induce a 1-1 correspondence between the limiting tangent hyperplanes to $X$ at the origin, and the limiting tangent hyperplanes to $X(0)\subset \IC^n$ at the origin? An affirmative answer means we can reduce the study of the limiting tangent hyperplanes to $X$ to a study of the limiting tangent hyperplanes of the fibers of the family.
 
 If $X$ is an ICIS, then there is a numerical criterion for $H$ not to be a limiting tangent hyperplane, namely that $e(JM(f)_H)=e(JM(f))$ (Cf. prop. 2.6 \cite{G0}).

\thm 6 Suppose $X$ is a family of ICIS, $(X_0,Y)$ satisfies WA, and suppose $e(JM(X(y)))$ independent of $Y$. Then $\iota$ induces a 1-1 correspondence between the limiting tangent hyperplanes to $X$ at $0$ and the limiting tangent hyperplanes in $\IC^n$ to $X\cap (0\times \IC^n)$.

\pf  Suppose $H$ is not a limiting tangent hyperplane to $X$, $H\supset Y$. A careful examination of the proof of Theorem 2.2 shows that $JM_z(f)_H$ is still a reduction of $JM(f)$, since the WA condition implies that $JM_y(f)\subset JM(f)^\dagger$. Restricting to $0\times \IC^n$, shows that $\iota(H)$ is not a limiting tangent hyperplane to $X(0)$, because $JM_z(f)_H$ remains a reduction of $JM_z(f)$.

Suppose $H\cap\IC^n$ is not a limiting tangent hyperplane to $X(0)$ at the origin. Then in 
$\O^p_{X(0),0}$,  $JM(f_0)_{\iota(H)}$ is a reduction of $JM(f_0)$, hence the multiplicity of the two modules are the same. The multiplicity of the family of modules $JM(f_y)_{\iota(H)}$ must be constant, because $e(M(f_y)_{\iota(H)},0)\ge e(JM(f_y),0)$, $e(M(f_y)_{\iota(H)},0)$ is upper semicontinuous and 
$e(JM(f_y),0)$ is constant by hypothesis. This means that $JM(f_y)_{\iota(H)}$ is a reduction of $JM(f_y)$ at the origin for all $y$, hence $\iota(H)$ is not a limiting tangent hyperplane to $X(y)$ at the origin.

Thus we have that the multiplicity of $JM(f_y)_{\iota(H)}$ at the origin is independent of $y$, and $JM_z(f)$ is fiberwise integrally dependent on $JM_z(f)_H$. The principle of specialization of integral dependence for modules (PSID) (\cite{GK}) then implies that $JM_z(f)_H$ is a reduction of  $JM_z(f)$ in 
$\O_{X,0}$. Whitney A then implies $JM_z(f)_H$ is a reduction of $JM(f)$, hence $H$ is not a limiting tangent hyperplane. Since $\iota$ induces a 1-1 correspondence between hyperplanes which are not tangent hyperplanes, it then follows that it induces a 1-1 correspondence between limiting tangent hyperplanes, which finishes the proof.

The condition of the theorem that $e(JM(X(y)))$ independent of $Y$ has a geometric interpretation--it means that there is no relative polar variety of $X$ of dimension equal to $Y$.

With Theorem 2.6, we can improve 2.2, 2.3 and 2.5. 

\cor 7 Suppose $X$ is a family of ICIS, $(X-Y,Y)$ satisfies W, then the following are equivalent 

1) $\{H_i\}$ is a W-generic sequence.

2)Each $\iota(H_i)$ is not a limiting tangent hyperplane at the origin of $X_{i-1}(0)$.

3) $e(JM(f_0),0)=e(JM(f)_{\iota(H)},0)$.

\pf If $X$ is a family of ICIS, with $X-Y$ smooth, $X-Y,Y$ Whitney, it follows that the multiplicity of $JM(f_y)$ is independent of $y$ at the origin.  The proof then proceeds by induction on $i$. By Theorem 2.6, the equivalent conditions 2) and 3)  imply that $H_i$ is not a limiting tangent hyperplane of $X_{i-1}$. Theorem 2.2 implies $H_i$ is generic for $X_{i-1}$, while  $H_i$ is not a limiting tangent hyperplane of $X_{i-1}$ also implies $S(X_i)\subset Y$. 

If we assume 1), then Theorem 2.2 implies that $H$ is not a limiting tangent hyperplane to $X$ at the origin, and 2.6 implies 2) and 3).

In \cite{G3} there is a numerical criterion for $H$ not to be a limit of tangent hyperplanes which holds for equidimensional spaces. It is based on the author's extension of the Buchsbaum-Rim multiplicity to modules of non-finite colength. Using this criterion and the generalization of the PSID in \cite{G4}, it is reasonable to expect that the analogue of 2.7 holds in general.

\sct References

\references
G0
T.  Gaffney, {\it Aureoles and integral closure of modules,} Stratifications, singularities and differential equations, II (Marseille, 1990; Honolulu, HI, 1990), 55--62, Travaux en Cours, 55, Hermann, Paris, 1997.

G1  
 T. Gaffney,
 {\it Integral closure of modules and Whitney equisingularity,}
 \invent 107 1992 301--22

G2
T. Gaffney
{\it Equisingularity of plane sections, $t\sb 1$ condition
 and the integral closure of modules}, Real and complex singularities
 (Sao Carlos, 1994), Pitman Res. Notes Math. Ser., Longman,
Harlow, vol. 333, 1995, p 95--111

G3
T.  Gaffney, {\it Generalized Buchsbaum-Rim multiplcities and a theorem of Rees}, Special issue in honor of Steven L. Kleiman. Comm. Algebra 31 (2003), no. 8, 3811--3827.

G4
	T. Gaffney, {\it The Multiplicity-Polar Formula and Equisingularity}, 
in preparation

GK
 T. Gaffney and S. Kleiman.
 {\it Specialization of integral dependence for modules}.
 Inventiones  math. 137, 541-574 1999

GK2  
T.  Gaffney and  S. L. Kleiman, {\it \Wf and integral dependence}, Real and complex singularities (S‹o Carlos, 1998), 33--45, Chapman and Hall/CRC Res. Notes Math., 412, Chapman and Hall/CRC, Boca Raton, FL, 2000.

KT1	
 S. Kleiman and A. Thorup,
  {\it A geometric theory of the Buchsbaum--Rim multiplicity,}
 \ja 167 1994 168--231
 
 LT
 D.T. L\^{e} and B. Teissier, {\it Limites d'espaces tangents en g\'{e}om\'{e}trie analytique}, Comm. Math. Helv., vol 63, 1988, p. 540--578

\endreferences

\bye